\input amstex
\documentstyle{amsppt}
\magnification=\magstep1                        
\hsize6.5truein\vsize8.9truein                  
\NoRunningHeads
\loadeusm

\magnification=\magstep1                        
\hsize6.5truein\vsize8.9truein                  
\NoRunningHeads
\loadeusm

\document
\topmatter

\title
On the oscillation of the modulus of the Rudin-Shapiro polynomials on the unit circle\\
\endtitle

\rightheadtext{the oscillation of the modulus of the Rudin-Shapiro polynomials on the unit circle}

\author Tam\'as Erd\'elyi
\endauthor

\address Department of Mathematics, Texas A\&M University,
College Station, Texas 77843, College Station, Texas 77843 \endaddress

\endthanks


\date February 28, 2018 
\enddate

\email terdelyi\@math.tamu.edu
\endemail

\abstract
In signal processing the Rudin-Shapiro polynomials have good autocorrelation properties
and their values on the unit circle are small. Binary sequences with low autocorrelation
coefficients are of interest in radar, sonar, and communication systems.
In this paper we study the oscillation of the modulus of the Rudin-Shapiro polynomials on
the unit circle. We also show that the Rudin-Shapiro polynomials $P_k$ and $Q_k$ of degree 
$n-1$ with $n := 2^k$ have $o(n)$ zeros on the unit circle. This should be compared with a 
result of B. Conrey, A. Granville, B. Poonen, and K. Soundararajan stating that for odd 
primes $p$ the Fekete polynomials $f_p$ of degree $p-1$ have asymptotically $\kappa_0 p$ 
zeros  on the unit circle,  where $0.500813>\kappa_0>0.500668$. Our approach is based heavily 
on the Saffari and Montgomery conjectures proved recently by B. Rodgers.  
We also prove that there are absolute constants $c_1 > 0$ and $c_2 > 0$ such that the $k$-th 
Rudin-Shapiro polynomials $P_k$ and $Q_k$ of degree $n-1 = 2^k-1$ have at least $c_2n$ zeros 
in the annulus 
$$\left \{z \in {\Bbb C}: 1 - \frac{c_1}{n} < |z| < 1 + \frac{c_1}{n} \right \}\,.$$

\endabstract

\endtopmatter

\head 1. Introduction and Notation \endhead

Let $D := \{z \in {\Bbb C}: |z| < 1\}$ denote the open unit disk of the complex plane.
Let $\partial D :=  \{z \in {\Bbb C}: |z| = 1\}$ denote the unit circle of the complex plane.
The Mahler measure $M_{0}(f)$ is defined for bounded measurable functions $f$ on $\partial D$ by 
$$M_{0}(f) := \exp\left(\frac{1}{2\pi} \int_{0}^{2\pi}{\log|f(e^{it})|\,dt} \right)\,.$$
It is well known, see [HL-52], for instance, that
$$M_{0}(f) = \lim_{q \rightarrow 0+}{M_{q}(f)}\,,$$
where
$$M_{q}(f) := \left( \frac{1}{2\pi} \int_{0}^{2\pi}{\left| f(e^{it}) \right|^q\,dt} \right)^{1/q}\,, 
\qquad q > 0\,.$$
It is also well known that for a function $f$ continuous on $\partial D$ we have 
$$M_{\infty}(f) := \max_{t \in [0,2\pi]}{|f(e^{it})|} = \lim_{q \rightarrow \infty}{M_{q}(f)}\,.$$
It is a simple consequence of the Jensen formula that
$$M_0(f) = |c| \prod_{j=1}^n{\max\{1,|z_j|\}}$$
for every polynomial of the form
$$f(z) = c\prod_{j=1}^n{(z-z_j)}\,, \qquad c,z_j \in {\Bbb C}\,.$$
See [BE-95, p. 271] or [B-02, p. 3], for instance. It will be convenient for us to introduce the notation 
$$M_{q}(S) := M_{q}(f)\,, \qquad 0 \leq q \leq \infty\,,$$
for functions $S$ defined on the period $K := {\Bbb R} \enskip(\text {mod}\,\, 2\pi)$ by 
$S(t) := f(e^{it})$, where $f$ is a bounded measurable functions $f$ on $\partial D$. 

Let ${\Cal P}_n^c$ be the set of all algebraic polynomials of degree at most $n$ with complex coefficients. 
Let ${\Cal T}_n$ be the set of all real (that is, real-valued on the real line) trigonometric polynomials 
of degree at most $n$. Finding polynomials with suitably restricted coefficients and maximal Mahler measure has 
interested many authors. The classes
$${\Cal L}_n := \left\{ f: \enskip f(z) = \sum_{j=0}^{n}{a_jz^j}\,, \quad a_j \in \{-1,1\} \right\}$$
of Littlewood polynomials and the classes
$${\Cal K}_n := \left\{ f: \enskip f(z) = \sum_{j=0}^{n}{a_jz^j}\,, \quad a_j \in {\Bbb C}, \enskip |a_j| =1 \right\}$$
of unimodular polynomials are two of the most important classes considered.
Observe that ${\Cal L}_n \subset {\Cal K}_n$ and
$$M_0(f) \leq M_2(f) = \sqrt{n+1}$$
for every $f \in {\Cal K}_n$.
Beller and Newman [BN-73] constructed unimodular polynomials $f_n \in {\Cal K}_n$  whose
Mahler measure $M_0(f_n)$ is at least $\sqrt{n}-c/\log n$.

Section 4 of [B-02] is devoted to the study of Rudin-Shapiro polynomials. 
Littlewood asked if there were polynomials $f_{n_k} \in {\Cal L}_{n_k}$ satisfying  
$$c_1 \sqrt{n_k+1}  \leq |f_{n_k}(z)| \leq c_2 \sqrt{n_k+1}\,, \qquad z \in \partial D\,,$$
with some absolute constants $c_1 > 0$ and $c_2 > 0$, see [B-02, p. 27] for a reference 
to this problem of Littlewood.
To satisfy just the lower bound, by itself, seems very hard, and no such sequence $(f_{n_k})$  
of Littlewood polynomials $f_{n_k} \in {\Cal L}_{n_k}$ is known. A sequence of Littlewood polynomials 
that satisfies just the upper bound is given by the Rudin-Shapiro polynomials. The Rudin-Shapiro 
polynomials appear in Harold Shapiro's 1951 thesis [S-51] at MIT and are sometimes called just 
Shapiro polynomials. They also arise independently in Golay's paper [G-51]. They are 
remarkably simple to construct and are a rich source of counterexamples to possible 
conjectures.

The Rudin-Shapiro polynomials are defined recursively as follows:
$$\split P_0(z) & :=1\,, \qquad Q_0(z) := 1\,, \cr 
P_{k+1}(z) & := P_k(z) + z^{2^k}Q_k(z)\,, \cr
Q_{k+1}(z) & := P_k(z) - z^{2^k}Q_k(z)\,, \cr \endsplit$$
for $k=0,1,2,\ldots\,.$ Note that both $P_k$ and $Q_k$ are polynomials of degree $n-1$ with $n := 2^k$
having each of their coefficients in $\{-1,1\}$.
In signal processing, the Rudin-Shapiro polynomials have good autocorrelation properties 
and their values on the unit circle are small. Binary sequences with low autocorrelation 
coefficients are of interest in radar, sonar, and communication systems.

It is well known and easy to check by using the parallelogram law that
$$|P_{k+1}(z)|^2 + |Q_{k+1}(z)|^2 = 2(|P_k(z)|^2 + |Q_k(z)|^2)\,, \qquad z \in \partial D\,.$$
Hence
$$|P_k(z)|^2 + |Q_k(z)|^2 = 2^{k+1} = 2n\,, \qquad z \in \partial D\,. \tag 1.1$$
It is also well known (see Section 4 of [B-02], for instance), that
$$Q_k(-z) = P_k^*(z) = z^{n-1}(P_k(1/z)\,, \qquad z \in \partial D\,,$$
and hence
$$|Q_k(-z)| = |P_k(z)|\,, \qquad z \in \partial D\,.i \tag 1.2$$
P. Borwein's book [B-02] presents a few more basic results on the Rudin-Shapiro
polynomials. Various properties of the Rudin-Shapiro polynomials are discussed in [B-73] by Brillhart 
and in [BL-76] by Brillhart, Lemont, and Morton.  
Obviously $M_2(P_k) = 2^{k/2}$ by the Parseval formula. In 1968 Littlewood
[L-68] evaluated $M_4(P_k)$ and found that $M_4(P_k) \sim (4^{k+1}/3)^{1/4}$.
The $M_4$ norm of Rudin-Shapiro like polynomials on $\partial D$ are studied in [BM-00]. 

P. Borwein and Lockhart [BL-01] investigated the asymptotic behavior of the mean
value of normalized $M_q$ norms of Littlewood polynomials for arbitrary $q > 0$.
They proved that
$$\lim_{n \rightarrow \infty} {\frac{1}{2^{n+1}} \, \sum_{f \in {\Cal L}_n}
{\frac{(M_q(f))^q}{n^{q/2}}}}= \Gamma \left( 1+ \frac q2 \right)\,.$$
In [C-15c] we proved that
$$\lim_{n \rightarrow \infty} {\frac{1}{2^{n+1}} \sum_{f \in {\Cal L}_n}
{\frac{M_q(f)}{n^{1/2}}}} = \left( \Gamma \left( 1+ \frac q2 \right) \right)^{1/q}$$
for every $q > 0$. In [CE-15c] we also proved the following result on the average Mahler measure of
Littlewood polynomials. We have
$$\lim_{n \rightarrow \infty} {\frac{1}{2^{n+1}} \sum_{f \in {\Cal L}_n}
{\frac{M_0(f)}{n^{1/2}}}} = e^{-\gamma/2}\,,$$
where
$$\gamma := \lim_{n \rightarrow \infty}{\left( \sum_{k=1}^n{\frac 1k - \log n} \right)} = 0.577215 \ldots$$
is the Euler constant and $e^{-\gamma/2} = 0.749306\ldots$.
These are analogues of the results proved earlier by Choi and Mossinghoff
[CM-11] for polynomials in ${\Cal K}_n$.
Let $K := {\Bbb R} \enskip(\text {mod}\,\, 2\pi)$.
Let $m(A)$ denote the one-dimensional Lebesgue measure of $A \subset K$.
In 1980 Saffari conjectured the following.

\proclaim{Conjecture 1.1}
Let $P_k$ and $Q_k$ be the Rudin-Shapiro polynomials of degree $n-1$ with $n := 2^k$. 
We have 
$$M_q(P_k) = M_q(Q_k) \sim \frac{2^{k+1)/2}}{(q/2+1)^{1/q}}$$
for all real exponents $q > 0$. Equivalently, we have 
$$\split & \lim_{k \rightarrow \infty} 
m{\left(\left\{t \in K: \left| \frac{P_k(e^{it})}{\sqrt{2^{k+1}}} \right|^2 \in [\alpha,\beta] \right\}\right)} \cr 
= \, & \lim_{k \rightarrow \infty}
m{\left(\left\{t \in K: \left| \frac{Q_k(e^{it})}{\sqrt{2^{k+1}}} \right|^2 \in [\alpha,\beta] \right\}\right)} 
= 2\pi(\beta - \alpha) \cr \endsplit$$
whenever $0 \leq \alpha < \beta \leq 1$. 
\endproclaim

This conjecture was proved for all even values of $q \leq 52$ by Doche [D-05]
and Doche and Habsieger [DH-04]. Recently B. Rodgers [R-16] proved Saffari's Conjecture 1.1 
for all $q > 0$. See also [EZ-17]. An extension of Saffari's conjecture is Montgomery's conjecture below.  

\proclaim{Conjecture 1.2}
Let $P_k$ and $Q_k$ be the Rudin-Shapiro polynomials of degree $n-1$ with $n := 2^k$.
We have
$$\split & \lim_{k \rightarrow \infty} 
m{\left(\left\{t \in K: \frac{P_k(e^{it})}{\sqrt{2^{k+1}}} \in E \right\}\right)} \cr
= \, & \lim_{k \rightarrow \infty}
m{\left(\left\{t \in K: \frac{Q_k(e^{it})}{\sqrt{2^{k+1}}} \in E \right\}\right)} 
= 2m(E) \cr \endsplit$$
for any measurable set $E \subset D := \{z \in {\Bbb C}: |z| < 1\}\,.$
\endproclaim

B. Rodgers [R-16] proved Montgomery's Conjecture 1.2 as well. 

Despite the simplicity of their definition not much is known
about the Rudin-Shapiro polynomials. It has been shown in [E-16c] fairly recently that 
the Mahler measure ($M_0$ norm) and the $M_\infty$ norm of the Rudin-Shapiro polynomials 
$P_k$ and $Q_k$ of degree $n-1$ with $n := 2^k$ on the unit circle of the complex plane 
have the same size, that is, the Mahler measure of the Rudin-Shapiro polynomials   
of degree $n-1$ with $n := 2^k$ is bounded from below by $cn^{1/2}$, where $c > 0$ is 
an absolute constant.

It is shown in this paper that the Rudin-Shapiro polynomials $P_k$ and $Q_k$ of degree $n-1$ 
with $n := 2^k$ have $o(n)$ zeros on the unit circle. We also prove that there are 
absolute constants $c_1 > 0$ and $c_2 > 0$ such that the $k$-th Rudin-Shapiro polynomials
$P_k$ and $Q_k$ of degree $n-1 = 2^k-1$ have at least $c_2n$ zeros in the annulus
$$\left \{z \in {\Bbb C}: 1 - \frac{c_1}{n} < |z| < 1 + \frac{c_1}{n} \right \}\,,$$
while there is an absolute constant $c > 0$ such that each of the functions
$\text {\rm Re}(P_k)$, $\text {\rm Re}(Q_k)$, $\text {\rm Im}(P_k)$, and $\text {\rm Im}(Q_k)$ 
has at least $cn$ zeros on the unit circle. The oscillation of $R_k(t) := |P_k(e^{it})|^2$ and 
$R_k(t) := |Q_k(e^{it})|^2$ on the period $[0,2\pi)$ is also studied. 
 
For a prime number $p$ the $p$-th Fekete polynomial is defined as
$$f_p(z) := \sum_{j=1}^{p-1}{\left( \frac jp \right)z^j}\,,$$
where
$$\left( \frac jp \right) =
\cases
1, \quad \text{if \enskip} x^2 \equiv j \enskip (\text {mod\,}p) \enskip
\text{has a nonzero solution,}
\\
0, \quad \text{if \enskip} p \enskip \text{divides} \enskip j\,,
\\
-1, \quad \text{otherwise}
\endcases$$
is the usual Legendre symbol. Since $f_p$ has constant coefficient $0$, it is not
a Littlewood polynomial, but $g_p$ defined by $g_p(z) := f_p(z)/z$ is a Littlewood
polynomial of degree $p-2$. Fekete polynomials are examined in detail in 
[B-02], [CG-00], [E-11], [E-12], [E-17], [EL-07], and [M-80]. In [CE-15a] and [CE-15b] the authors examined 
the maximal size of the Mahler measure and the $L_p$ norms of sums of $n$ monomials on the unit circle 
as well as on subarcs of the unit circles. In the constructions appearing in [CE-15a] 
properties of the Fekete polynomials $f_p$ turned out to be quite useful.
In [CG-00] B. Conrey, A. Granville, B. Poonen, and K. Soundararajan proved that
for an odd prime $p$ the Fekete polynomial $f_p(z)=\sum^{p-1}_{j=0} \big({j \over p}\big) z^j$
(the coefficients are Legendre symbols) has $\sim \kappa_0 p$  zeros on the unit circle, where
$0.500813>\kappa_0>0.500668$. So Fekete polynomials are far from having only $o(p)$ zeros
on the unit circle.

Mercer [M-06a] proved that if a Littlewood polynomial $f \in {\Cal L}_n$ of the form $f(z) = \sum_{j=0}^n{a_jz^j}$
is skew-reciprocal, that is, $a_j = (-1)^ja_{n-j}$ for each $j=0,1,\ldots,n$, then it has
no zeros on the unit circle. However, by using different elementary methods it was observed
in both [E-01] and [M-06a] that if a Littlewood polynomial $P$ of the form (1.1) is self-reciprocal,
that is, $a_{j,n} = a_{n-j}$ for each $j=0,1,\ldots,n$, $n \geq 1$, then it has at least one zero
on the unit circle.

There are many other papers on the zeros of constrained polynomials. Some of them are 
[BP-32], [BE-97], [BE-01], [BE-07], [BE-08a], [BE-08b], [BE-99], [BE-13], [B-97], [D-08], [E-08a], [E-08b], [E-16a], [E-16b],    
[L-61], [L-64], [L-66a], [L-66b], [L-68], [M-06b], [Sch-32], [Sch-33], [Sz-34], and [TV-07]. 

\head 2. New Results \endhead

Let $P_k$ and $Q_k$ be the Rudin-Shapiro polynomials of degree $n-1$ with $n := 2^k$.
Let either $R_k(t) := |P_k(e^{it})|^2$ or $R_k(t) := |Q_k(e^{it})|^2$. Let $\gamma := \sin^2(\pi/8)$.
We use the notation 
$$\|g\|_A := \sup_{x \in A}{|g(x)|}$$
for a complex-valued function $g$ defined on a set $A \subset {\Bbb R}$.

\proclaim{Theorem 2.1}
$P_k$ and $Q_k$ have $o(n)$ zeros on the unit circle.
\endproclaim

The proof of Theorem 2.1 will follow by combining the recently proved Saffari's conjecture
stated as Conjecture 1.1 and the theorem below. Let $K := {\Bbb R} \enskip(\text {mod}\,\, 2\pi)$.

\proclaim{Theorem 2.2} If $S \in {\Cal T}_n$ is of the form $S(t) = |f(e^{it})|^2$,
where $f \in {\Cal P}_n^c$, and $f$ has at least $k$ zeros (counted with multiplicities) in $K$, then
$$m(\{t \in K: |S(t)| \leq \alpha \|S\|_K\}) \geq \frac{\sqrt{\alpha}}{e} \, \frac kn$$
for every $\alpha \in (0,1)$, where $m(A)$ denotes the one-dimensional Lebesgue measure of a measurable set 
$A \subset K$.
\endproclaim

\proclaim{Theorem 2.3}
There is an absolute constant $c > 0$ such that each of the functions
$\text {\rm Re}(P_k)$, $\text {\rm Re}(Q_k)$, $\text {\rm Im}(P_k)$, and $\text {\rm Im}(Q_k)$ has
at least $cn$ zeros on the unit circle for every $n=2^k-1 \geq 1$.
\endproclaim

\proclaim{Theorem 2.4}
There is an absolute constant $c > 0$ such that the equation $R_k(t) = \eta n$ has at most 
$c\eta^{1/2}n$ solutions (counted with multiplicities) in $K$  for every $\eta \in (0,1]$ and sufficiently large 
$k \geq k_\eta$,  while the equation $R_k(t) = \eta n$  has at most $c(2-\eta)^{1/2}n$ solutions 
(counted with multiplicities) in $K$ for every $\eta \in [1,2)$ and sufficiently large $k \geq k_\eta$.
\endproclaim

\proclaim{Theorem 2.5}
The equation $R_k(t) = \eta n$ has at least $(1-\varepsilon)\eta n/2$ distinct solutions in $K$
for every $\eta \in (0,2\gamma)$, $\varepsilon > 0$, and sufficiently large $k \geq k_{\eta,\varepsilon}$,  
while the equation $R_k(t) = \eta n$ has at least $(1-\varepsilon)(2-\eta) n/2$ distinct solutions in $K$ 
for every $\eta \in (2-\gamma,2)$,  $\varepsilon > 0$, and sufficiently large $k \geq k_{\eta,\varepsilon}$.
\endproclaim

\proclaim{Theorem 2.6}
There s an absolute constants $A>0$ such that the equation $R_k(t) = (1+\eta)n$ has at least $An^{0.36}$ 
distinct solutions in $K$ whenever $\eta$ is real and $|\eta| < 2^{-11}$.
\endproclaim

\proclaim{Theorem 2.7}
There are absolute constants $c_1 > 0$ and $c_2 > 0$ such that $P_k$ and $Q_k$ have at least 
$c_2n$ zeros in the annulus
$$\left \{z \in {\Bbb C}: 1 - \frac{c_1}{n} < |z| < 1 + \frac{c_1}{n} \right \}\,.$$
\endproclaim

We note that for every $c \in (0,1)$ there is an absolute constant $c_3 > 0$ depending only on $c$ such that
every $U_n \in {\Cal P}_n^c$ of the form
$$U_n(z) = \sum_{j=0}^n{a_jz^j}\,, \qquad |a_0| = |a_n| = 1\,, \quad a_j \in {\Bbb C}\,, \quad |a_j| \leq 1\,,$$
has at least $cn$ zeros in the annulus
$$\left \{z \in {\Bbb C}: 1-\frac{c_3 \log n}{n} < |z| < 1+\frac{c_3 \log n}{n} \right\}\,. \tag 2.1$$
See Theorem 2.1 in [E-01].
  
On the other hand, there is an absolute constant $c_3 > 0$ such that
for every $n \in {\Bbb N}$ there is a polynomial $U_n \in {\Cal K}_n$ having no zeros in the annulus (2.1).
See Theorem 2.3 in [E-01]. So in Theorem 2.7 some special properties, in addition to being Littlewood polynomials, 
of the Rudin-Shapiro polynomials must be exploited.  

A key to the proof of Theorem 2.7 is the result below.

\proclaim{Theorem 2.8}
Let $t_0 \in {\Bbb R}$. There is an absolute constant $c_4 > 0$ depending only on $c > 0$ such that
$P_k$ has at least one zero in the disk
$$\left \{z \in {\Bbb C}: |z - e^{it_0}| < \frac{c_4}{n} \right\}\,,$$
whenever
$$R_k^{\prime}(t_0) \geq cn^2\,, \qquad R_k(t) = P_k(e^{it})P_k(e^{-it})\,.$$
\endproclaim

\proclaim{Problem 2.9}
Is there an absolute constant $c > 0$ such that the equation $R_k(t) = \eta n$ has at least
$c\eta n$ distinct solutions in $K$ for every $\eta \in (0,1)$ and sufficiently large $k \geq n_{\eta}$?
In other words, can Theorem 2.5 be extended to all $\eta \in (0,1)$?
\endproclaim

We note that it follows from $Q_k(z) = P_k^*(-z)$, $z \in {\Bbb C}$, that the products $P_kQ_k$ have
at least $n-1$ zeros in the closed unit disk and at least $n-1$ zeros outside the open unit disk.      
So in the light of Theorem 2.1 the products $P_kQ_k$ have asymptotically $n$ zeros in the open
unit disk. However, as far as we know, the following questions are open. 

\proclaim{Problem 2.10}
Is there an absolute constant $c > 0$ such that $P_k$ has at least $cn$ zeros in the open unit disk?
\endproclaim

\proclaim{Problem 2.11}
Is there an absolute constant $c > 0$ such that $Q_k$ has at least $cn$ zeros in the open unit disk?
\endproclaim

\proclaim{Problem 2.12}
Is it true that both $P_k$ and $Q_k$ have asymptotically half of their zeros in the open unit disk?
\endproclaim

\proclaim{Problem 2.13}
Is it true that if $n$ is odd then $P_k$ has a zero on the unit circle ${\partial D}$ only at $-1$ and
$Q_k$ has a zero on the unit circle ${\partial D}$ only at $1$, while if $n$ is even then neither $P_k$ 
nor $Q_k$ has a zero on the unit circle?
\endproclaim

As for $k \geq 1$ both $P_k$ and $Q_k$ have odd degree, both $P_k$ and $Q_k$ have at least one  
real zero. The fact that for $k \geq 1$ both $P_k$ and $Q_k$ have exactly one real zero was proved 
in [B-73].

\head 3. Lemmas \endhead 

Let, as before, $P_k$ and $Q_k$ be the Rudin-Shapiro polynomials of degree $n-1$ with $n := 2^k$.
To prove Theorem 2.1 we need the lemma below that is proved in [BE-95, E.11 of Section 5.1
on pages 236--237].

\proclaim{Lemma 3.1} If $S \in {\Cal T}_n$, $t_0 \in K$, and $r > 0$, then $S$ has
at most $enr|S(t_0)|^{-1} \|S\|_K$ zeros in the interval $[t_0 - r, t_0 + r]$.
\endproclaim

Our next lemma is in [E-18], and for the sake of completeness we present its short 
proof in  

Our next lemma is stated as Lemma 3.5 in [E-16c], where its proof may also be found.

\proclaim{Lemma 3.2}
If $\gamma := \sin^2(\pi/8)$ and
$$z_j := e^{it_j}\,, \quad t_j := \frac{2\pi j}{n}\,, \qquad j \in {\Bbb Z}\,,$$
then
$$\max \{|P_k(z_j)|^2,|P_k(z_{j+r})|^2\} \geq \gamma 2^{k+1} = 2\gamma n\,, \quad r \in \{-1,1\}\,,$$
for every $j=2u$, $u \in {\Bbb Z}$.
\endproclaim

By Lemma 3.2, for every $n = 2^k$  there are
$$0 \leq \tau_1 < \tau_2 < \cdots < \tau_m < \tau_{m+1} := \tau_1 + 2\pi$$
such that
$$\tau_j - \tau_{j-1} = \frac{2\pi l}{n}\,, \qquad l \in \{1,2\}\,,$$
and with
$$a_j := e^{i\tau_j}, \qquad j=1,2,\ldots,m+1, \tag 3.4$$
we have
$$|P_k(a_j)|^2 \geq 2\gamma n\,, \qquad j=1,2,\ldots,m+1\,.$$
(Moreover, each $a_j$ is an $n$-th root of unity.)  

Our next lemma is stated and proved as Lemma 3.4 in [E-18].

\proclaim{Lemma 3.3}
There is an absolute constant $c_3 > 0$ such that
$$\mu := \left|\left \{j \in \{2,3,\ldots,m+1\}: \min_{t \in [\tau_{j-1},\tau_j]}{R_k(t)} \leq \varepsilon \right \}\right| 
\leq c_3n\varepsilon^{1/2}$$
for every sufficiently large $n=2^k \geq n_\varepsilon$, $k=1,2,\ldots$, and $\varepsilon > 0$.
\endproclaim

Our next lemma is based on the work of M. Taghavi [T-97] and gives an upper bound for the 
so-called autocorrelation coefficients of the Rudin-Shapiro polynomials. 

\proclaim{Lemma 3.4}
If  
$$|P_k(e^{it})|^2 = \sum_{j=-n+1}^{n-1}{a_jz^j}$$
($a_0=n$, $a_j = a_{-j}$, $j \geq 1$), then 
$$\max_{1 \leq j \leq n-1}{|a_j|} \leq Cn^{0.8190}$$
with an absolute constant $C > 0$.
\endproclaim

In fact, Taghavi [T-97] claimed 
$$\max_{1 \leq j \leq n-1}{|a_j|} \leq (3.2134)n^{0.7303}\,.$$
However, as Allouche and Saffari observed, in his proof Taghavi used an incorrect 
statement saying that the spectral radius of the product of some matrices is independent 
of the order of the factors. So what he ended up with cannot be viewed as a correctly 
proved result. Building on what is correct in [T-97] Stephen Choi made some computations  
leading to the above correct form of Taghavi's upper bound on the autocorrelation 
coefficients of the Rudin-Shapiro polynomials. The correction based on Choi's 
computations will be the subject of a forthcoming note [AC-17] perhaps even in a more 
optimized form.

Using the notation of Lemma 3.4 Taghavi [T-96] claims also that
$$\max_{1 \leq j \leq n-1}{|a_j|} \geq Cn^{0.73}$$ for all $n := 2^k$ with an 
absolute constant $C>0$. 

Our next lemma is due to Littlewood, see [Theorem 1 in L-66a].

\proclaim{Lemma 3.5}
If $S \in {\Cal T}_n$ of the form 
$$S(t) = \sum_{m=0}^n{b_m \cos(mt + \alpha_m)}\,, \qquad b_m, \alpha_m \in {\Bbb R}\,, \tag 3.1$$
satisfies 
$$M_1(S) \geq c\mu\,, \qquad \mu := M_2(S)\,,$$
where $c > 0$ is a constant, $b_0 = 0$,  
$$s_{\lfloor n/h \rfloor} = \sum_{m=1}^{\lfloor n/h \rfloor}{\frac{b_m^2}{\mu^2}} \leq 2^{-9}c^6$$
for some constant $h>0$, and $v \in {\Bbb R}$ satisfies 
$$|v| \leq V = 2^{-5}c^3\,,$$
then
$${\Cal N}(S,v) > Ah^{-1}c^5n\,,$$
where ${\Cal N}(S,v)$ denotes the number of real zeros of $S-v\mu$ in $(-\pi,\pi)$, and $A > 0$ is 
an absolute constant.
\endproclaim 

Our next lemma is a key to prove Theorem 2.7. It is an extension of Theorem 1 in [E-02] establishing 
the right Bernstein inequality for trigonometric polynomials $S \in {\Cal T}_n$ not vanishing in the strip 
$$\{z \in {\Bbb C}: |\text {\rm Im}(z)| < r\}\,, \qquad 0 < r \leq 1\,.$$ 

\proclaim{Lemma 3.6} Let $0 < r \leq 2$. We have
$$|S^{\prime}(a)| \leq 5e \sqrt{\frac{2n}{r}} \, \|S\|_K$$
for every $S \in {\Cal T}_n$ having no zeros in the disk $D(a,r) := \{z \in {\Bbb C}: |z-a| < r\}$.
\endproclaim

For the proof of Lemma 3.6 we need the lemma below.

\proclaim{Lemma 3.7}
Let $2/n < r \leq 2$. We have
$$|S(z)| \leq e\|S\|_K$$
for every $S \in {\Cal T}_n$ having no zeros in the disk $D(a,r)$ centered at $a$ of radius $r$,
and for every $z$ in the square
$$\{z = x+iy: x \in [a - \rho, a + \rho], y \in [-\rho,\rho]\}$$
with $\rho := \displaystyle{\frac 15 \left( \frac{r}{2n} \right)^{1/2}}$.
\endproclaim

For the proof of Lemma 3.7 we need the lemma below.

\proclaim{Lemma 3.8}
Let $1/n \leq r \leq 1$. We have
$$|S(x+iy)| \leq e\|S\|_K$$
for every $S \in {\Cal T}_n$ having no zeros in the disk $D(x,r)$, and for every $y \in [-\rho,\rho]$ with
$\rho := \displaystyle{\frac 15 \left( \frac{r}{n} \right)^{1/2}}$.
\endproclaim

For the proof of Lemma 3.8 we need the lemma below stated as Lemma 4.3 in [E-98].

\proclaim{Lemma 3.9} Let $0 < s \leq \lambda \leq 1$. We have
$$\|f\|_{[-1,1+s]} \leq \exp\big(8n\lambda^{-1/2}s\big)\,\|f\|_{[-1,1]}$$
for every $f \in {\Cal P}_n^c$ having no zeros in the disk $D(1-\lambda,\lambda)$ centered at $1-\lambda$ of radius $\lambda$.
\endproclaim

\demo{Proof of Lemma 3.8} It is sufficient to prove the lemma for $x=0$, since for $x \neq 0$
we can study the polynomial $\widetilde{S} \in {\Cal T}_n$ defined by $\widetilde{S}(t) := S(t-x)$
having no zeros in the disk $D(0,r)$. Associated with $S \in {\Cal T}_n$ we define $U \in {\Cal T}_{2n}$ 
by $U(t):=S(t)S(-t)$. Observe that $U$ is an even trigonometric polynomial of degree at most $2n$, hence 
we can define $f \in {\Cal P}_{2n}^c$ (in fact, with real coefficients) by
$$f(\cos t):= U(t), \qquad t \in {\Bbb C}\,.$$
Assume that $S$, and hence $U$, has no zeros in the disk $D(0,r)$.
We show that $f$ has no zeros in the disk $D(1,2\lambda)$ centered at $1$ of radius $2\lambda := r^2/4$.
Indeed, as $S$, and hence $U$, has no zeros in the disk $D(0,r)$, $f$ has no zeros in the 
region $H := \{u = \cos t: t \in D(0,r)\}$ bounded by the curve $\Gamma := \partial H := \{u = \cos t: |t| = r\}$.
As $\cos 0 = 1$, $\Gamma$ goes around $1$ at least once by the Argument Principle. 
Observe that if $z = \cos t \in \Gamma$, then $|t| = r \leq 1$ implies that
$$|1 - \cos t| =  \left| \sum_{k=1}^\infty{\frac{|t|^{2k}}{(2k!}} \right| 
\geq |t|^2 \left( \frac 12 - \sum_{k=1}^\infty{\frac{|t|^{2k}}{(2k+2)!}} \right) \geq \frac{|t|^2}{4} = \frac{r^2}{4}\,,$$
and hence $H$ contains the disk $D(1,2\lambda) = D(1,r^2/4)$. In conclusion, $f$ has no zeros in the disk $D(1,2\lambda)$ as 
we claimed. 

Using Lemma 3.9 with $\displaystyle{\lambda := \frac{r^2}{8}}$ and $s:=u-1$, we have
$$\split |f(u)| & \leq \exp(8(2n)\lambda^{-1/2}(u-1))\|f\|_{[-1,1]} \cr 
& \leq \exp(48nr^{-1}(u-1))\|f\|_{[-1,1]}\,, \qquad 0 < s := u-1 \leq \lambda = \frac{r^2}{8} \leq 1\,. \cr \endsplit \tag 3.1$$
Now let 
$$\rho := \frac 15 \left( \frac rn \right)^{1/2}, \qquad y \in [-\rho,\rho]\,, \qquad  u := \cosh y\,.$$ 
Then $1/n \leq r \leq 1$ and $\displaystyle{\lambda := \frac{r^2}{8}}$ imply that 
$$u-1 = \cosh y -1 \leq y^2 \leq \rho^2 = \frac{r}{25n} < \frac{r^2}{8} = \lambda\,.$$
Using (3.1) we have
$$\split |S(iy)|^2 & = |U(iy)| = |f(\cos(iy))| = |f(\cosh y)| \leq \exp(48nr^{-1}(\cosh y -1))\|f\|_{[-1,1]}  \cr  
& \leq \exp(48nr^{-1}\rho^2)\|f\|_{[-1,1]} \leq \exp \left(48nr^{-1}\frac{r}{25n} \right)\|f\|_{[-1,1]} \cr  
& \leq e^2\|S\|_K^2 \cr \endsplit$$
for every $y \in [-\rho,\rho]$ with $\rho := \displaystyle{\frac 15 \left( \frac rn \right)^{1/2}}$, $1/n \leq r \leq 1$.
\qed \enddemo

\demo{Proof of Lemma 3.7}
Observe that if $S \in {\Cal T}_n$ has no zeros in the disk $D(a,r)$ centered at $a$ of radius $r$,
then it has no zeros in the disks $D(b,r/2)$ centered at $b \in [a-r/2,a+r/2]$ of radius $r/2$.
Observe that $2/n \leq r \leq 2$ implies $1/n \leq r/2 \leq 1$. Using Lemma 3.8 we obtain that
$$|S(x+iy)| \leq e\|S\|_K$$
for every $S \in {\Cal T}_n$ having no zeros in the disk $D(a,r)$ centered at $a$ of radius $r$,
and for every
$$z \in \{z=x+iy: x \in [a-r/2,a+r/2], y \in [-\rho,\rho]\}$$
with $\rho := \displaystyle{\frac 15 \left( \frac{r}{2n} \right)^{1/2}}$. As $2/n \leq r \leq 2$ implies $0 < \rho < r/2$, the lemma follows.
\qed \enddemo

\demo{Proof of Lemma 3.6}
If $2/n \leq r \leq 2$ then using Cauchy's integral formula and Lemma 3.7, we obtain
$$\split |S^\prime(a)| & \leq \left| \frac{1}{2\pi} \int_{|\zeta - a| = \rho}{\frac{S(\zeta)}{(\zeta-a)^2} \,d\zeta}\right| 
\leq e\rho^{-1} \|S\|_K \cr 
& \leq 5e \left( \frac{2n}{r} \right)^{1/2} \, \|S\|_K \cr \endsplit$$
for every $S \in {\Cal T}_n$ having no zeros in the disk $D(a,r) := \{z \in {\Bbb C}: |z-a| < r\}$.
If $r< 2/n$ then the classical Bernstein inequality valid for all $S \in {\Cal T}_n$ gives the lemma.
\qed \enddemo

\head 4. Proofs of the Theorems \endhead

\demo{Proof of Theorem 2.2}
Let $S \in {\Cal T}_n$ be of the form $S(t) = |f(e^{it})|^2$, where $f \in {\Cal P}_n^c$. We define 
$U \in {\Cal T}_n$ and $V \in {\Cal T}_n$ by 
$$U(t) := \text {\rm Re}(f(e^{it})) \qquad \text {\rm and} \qquad V(t) := \text {\rm Im}(f(e^{it}))\,, \qquad t \in K\,.$$
Then
$$S(t) = |f(e^{it})|^2 = U(t)^2 + V(t)^2\,, \qquad t \in K\,. \tag 4.1$$
Suppose $S \in {\Cal T}_n$ defined by $S(t) = |f(e^{it})|^2$ has at least $u$ zeros in $K$, 
and let $\alpha \in (0,1)$. Then 
$$\{t \in K: |S(t)| \leq \alpha \|S\|_K\}$$
can be written as the union of pairwise disjoint intervals 
$I_j$, $j=1,2,\ldots,m$. Each of the intervals $I_j$ contains a point $y_j \in I_j$ such that 
$$|S(y_j)| = \alpha \|S\|_K\,.$$ 
Hence, (4.1) implies that for each $j = 1,2,\ldots,m$, we have either  
$$|U(y_j)| \geq \left(\frac{\alpha}{2}\right)^{1/2} \|f\|_K \geq \left(\frac{\alpha}{2}\right)^{1/2} \|U\|_K \tag 4.2$$
or
$$|V(y_j)| \geq \left(\frac{\alpha}{2}\right)^{1/2} \|f\|_K \geq \left(\frac{\alpha}{2}\right)^{1/2} \|V\|_K\,. \tag 4.3$$
Also, each zero of $S$ lying in $K$ is contained in one of the intervals $I_j$. 
Let $\mu_j$ denote the number of zeros of $S$ lying in $I_j$. Since $S \in {\Cal T}_n$ 
has at least $u$ zeros in $K$, so do $U \in {\Cal T}_n$ and $V \in {\Cal T}_n$, and we have 
$\sum_{j=1}^m{\mu_j} \geq u$. Note that Lemma 3.1 applied to $U \in {\Cal T}_n$ yields that 
$$\mu_j \leq en|I_j| \left(\left(\frac{\alpha}{2}\right)^{1/2}\|U\|_K \right)^{-1} \|U\|_K = 
\frac{e\sqrt{2}\,n}{\alpha^{1/2}} \, |I_j|$$
for each $j = 1,2,\ldots,m$ for which (4.2) holds. 
Also, Lemma 3.1 applied to $V \in {\Cal T}_n$ yields that 
$$\mu_j \leq en|I_j| \left(\left(\frac{\alpha}{2}\right)^{1/2}\|V\|_K \right)^{-1} \|V\|_K = 
\frac{e\sqrt{2}\,n}{\alpha^{1/2}} \, |I_j|$$
for each $j = 1,2,\ldots,m$ for which (4.3) holds. Hence
$$\mu_j \leq \frac{e\sqrt{2}\,n}{\alpha^{1/2}}|I_j|\,, \qquad j = 1,2,\dots,m\,.$$
Therefore
$$u \leq \sum_{j=1}^m{\mu_j} \leq \frac{e\sqrt{2}\,n}{\alpha^{1/2}} \sum_{j=1}^m{|I_j|} 
= \frac{e\sqrt{2}\,n}{\alpha^{1/2}} \, m(\{t \in K: |S(t)| \leq \alpha |S\|_K\})\,,$$
and the lemma follows. 
\qed \enddemo

\demo{Proof of Theorem 2.1}
We show that the $P_k$ has $o(n)$ zeros on the unit circle, where $n = 2^k-1$. The proof of the fact that $Q_k$ has $o(n)$ 
zeros on the unit circle is analogous. Suppose to the contrary that there are $\varepsilon > 0$ and an increasing sequence 
$(k_j)_{j=1}^{\infty}$ of positive integers such that the Rudin-Shapiro polynomials $P_{k_j}$  have at least $\varepsilon n_j$ 
zeros on the unit circle, where $n_j := 2^{k_j}$ for each $j=1,2,\ldots$. Then $P_{k_j}$ has at least one zero on the unit circle 
and hence (1.1) and (1.2) imply that   
$$\|P_{k_j}(e^{it})\|_K^2 = 2^{k_j+1}\,. \tag 4.4$$
Then Theorem 2.2 implies that  
$$m(\{t \in K: |P_{k_j}(t)|^2 \leq \alpha \|P_{k_j}\|_K^2\}) \geq \frac{\sqrt{\alpha}}{e} \, \frac{\varepsilon n_j}{n_j}
= \, \frac{\varepsilon \sqrt{\alpha}}{e}$$
for every $\alpha \in (0,1)$ and $j=1,2,\ldots$. Hence,
$$\liminf_{j \rightarrow \infty} m(\{t \in K: |P_{k_j}(e^{it})|^2 \leq \alpha \|P_{k_j}(e^{it})\|_K^2\}) 
\geq \frac{\varepsilon \sqrt{\alpha}}{e}  \tag 4.5$$
for every $\alpha \in (0,1)$.
On the other hand, Conjecture 1.1 proved in [R-16] combined with (4.4) imply that 
$$\lim_{j \rightarrow \infty}
m{\left(\left\{t \in K: |P_{k_j}(e^{it})|^2 \leq \alpha \|P_{k_j}(e^{it})\|_K^2\right\}\right)} = 2\pi \alpha \tag 4.6$$
for every $\alpha \in (0,1)$. Combining (4.5) and (4.6) we obtain
$$\frac{\varepsilon \sqrt{\alpha}}{e} \leq 2\pi \alpha\,,$$
that is, $\varepsilon/e \leq 2\pi \sqrt{\alpha}$ for every $\alpha \in (0,1)$, a contradiction. 
\qed \enddemo

\demo{Proof of Theorem 2.3}
We prove that there is an absolute constant $c > 0$ such that $\text {\rm Re}(P_k)$ has at least $cn$ zeros 
on the unit circle; the fact that each of the functions $\text {\rm Re}(Q_k)$, $\text {\rm Im}(P_k)$, and 
$\text {\rm Im}(Q_k)$ has at least $cn$ zeros on the unit circle can be proved similarly. 
Let, as before, $K := {\Bbb R} \enskip(\text {mod}\,\, 2\pi)$. Let 
$${\Cal A}_n := \left\{f: f(t) =\sum_{j=1}^n{\cos(jt + \alpha_j)}\,, \enskip \alpha_j \in {\Bbb R}\right \}\,.$$
Let $S \in {\Cal A}_{n-1}$ with $n := 2^k$ be defined by
$$S(t) := \text {\rm Re}(P_k(e^{it})) - 1\,.$$
We have
$$\mu = M_2(S) := \left( \frac{1}{2\pi} \int_K{|S(t)|^2 \, dt} \right)^{1/2} = \left( \frac {n-1}{2} \right)^{1/2}.$$ 
Let ${\Cal N}(S,v)$ be the number of real roots of $S - v\mu$ in $[-\pi,\pi)$.
Observe that (1.1) implies that $|S(t)| \leq (2n)^{1/2} + 1 \leq 2(n-1)^{1/2}$ for all $t \in K$ and 
$n = 2^k-1 \geq 3$, and hence
$$\split M_1(S) = & \, \frac{1}{2\pi} \int_{0}^{2\pi}{|S(t)|\,dt} \geq 
\frac{1}{2\pi} \frac{1}{2(n-1)^{1/2}} \int_{0}^{2\pi}{|S(t)|^2\,dt} \cr
= & \, \frac{1}{2(n-1)^{1/2}} \frac{n-1}{2} = \frac{(n-1)^{1/2}}{4} = \frac{\mu}{2\sqrt{2}}\,. \cr \endsplit$$ 
for all $n = 2^k-1 \geq 3$.
Thus, applying Lemma 3.5 with $h = 2^9c^{-6}$ we can deduce that there is an absolute constant $A>0$ such that   
$$S(t) + 1 = \text {\rm Re}(P_k(e^{it}))$$ 
has at least $A2^{-16} 2^{-33/2} n = A2^{-65/2} n$ zeros in $[-\pi,\pi)$ 
whenever 
$$2^{-5}c^3\mu = 2^{-5}c^3\sqrt{(n-1)/2} \geq 1\,.$$ 
This finishes the proof when $n = 2^k-1$ is sufficiently large. On the other hand, if $k \geq 1$ then 
$P_k$ always has  at least one zero in the closed unit disk, hence $\text {\rm Re}(P_k(e^{it}))$ has at least 
two zeros in $[-\pi,\pi)$.  
\qed \enddemo

\demo{Proof of Theorem 2.4}
The proof is a combination of Lemmas 3.1, 3.2, and 3.3. Recalling (1.2) we can observe that 
without loss of generality we may assume that $\eta \in (0,1]$, that is, it is sufficient to prove 
only the first statement of the theorem. As the trigonometric polynomial $R_k(t) - \eta n$ of degree $n-1$ 
has at most $2(n-1)$ zeros in $K$,  without loss of generality we may assume 
also that $\eta < \gamma/2$, where $\gamma := \sin^2(\pi/8)$ as before. In the light of Lemma 3.3 it is sufficient 
to prove that there is an absolute constant $c > 0$ such that the equation $R_k(t) = \eta n$ has at most 
$c$ solutions in the interval $[\tau_{j-1},\tau_j]$ for every $j \in \{2,3,\ldots,m+1\}$ for which  
$$\min_{t \in [\tau_{j-1},\tau_j]}{R_k(t)} \leq \eta n\,.$$ 
However, this follows from Lemmas 3.1 combined with Lemma 3.2.  
\qed \enddemo

\demo{Proof of Theorem 2.5}
Recalling (1.1), without loss of generality we may assume that $\eta \in (0,2\gamma)$. Let
$$I_j := \left[ \frac{(2j-2)\pi}{n},\frac{2j\pi}{n} \right)\,, \qquad j=1,2,\ldots,n\,.$$
By Saffari's Conjecture 1.1 proved by Rodgers [R-16] we have
$$m(\{t \in K: R_k(t) \leq \eta n\}) > \pi (1-\varepsilon)\eta$$
for every $\eta \in (0,1)$, $\varepsilon > 0$, and sufficiently large $k \geq k_{\eta,\varepsilon}$.
Hence, with the notation 
$$A_\eta := \{t \in K: R_k(t) \leq \eta n\}\,,$$
there are at least $(1-\varepsilon)\eta n/2$ distinct values of $j \in \{1,2,\ldots,n\}$ such that 
$A_\eta \cap I_j \neq \emptyset$ for every $\eta \in (0,1)$ and sufficiently large $k \geq k_{\eta,\varepsilon}$. 
On the other hand, by Lemma 3.2, for each $j \in \{1,2,\ldots,n\}$ 
there is a $t_j \in I_j$ such that $R_k(t_j) \geq 2\gamma n$. Hence by the Intermediate 
Value Theorem there are at least $(1-\varepsilon)\eta n/2$ distinct values of $j \in \{1,2,\ldots,n\}$ 
for which there is a $\tau_j \in I_j$ such that $R_k(\tau_j) = \eta n$ 
for every $\eta \in (0,2\gamma)$, $\varepsilon > 0$, and sufficiently large $k \geq k_{\eta,\varepsilon}$. 
\qed \enddemo

\demo{Proof of Theorem 2.6}
Let $S_n \in {\Cal T}_{n-1}$ be defined by 
$$S_n(t) := R_k(t)-n = |P_k(e^{it})|^2-n = \sum_{j=-n+1}^{n-1}{a_jz^j}-n\,.$$ 
We show that $S := S_n$ satisfies the assumptions of Lemma 3.5 with $c=1/4$ and 
$h:=n^{0.64}$ if $n=2^k$ is sufficiently large.    
Clearly, $S_n$ is of the form (3.1) with $b_0=0$, $b_m = a_m/2$, and $\gamma_m=0$ for $m=1,2,\ldots,n-1$.
As it is already mentioned in Section 1, Littlewood [L-68] evaluated $M_4(P_k)$ and found that 
$M_4(P_k) \sim (4^{k+1}/3)^{1/4} = (4n^2/3)^{1/4}$. Hence $\mu := M_2(S_n) \sim (1/3)^{1/2}n$.   
Also, it follows from (1.1) that $M_{\infty}(S_n) \leq n$, hence  
$$(1/3)n^2 \sim (M_2(S_n))^2 \leq M_1(S_n)M_{\infty}(f) \leq nM_1(S_n)$$  
implies that $M_1(S_n) \geq (1/4)n$ if $n=2^k$ is sufficiently large. 
Now Lemma 3.4, $b_0=0$, $b_m = 2a_m \in \{-2,2\}$, $m=1,2,\ldots,n-1$, and $h:=n^{0.64}$ imply that
$$\split s_{\lfloor (n-1)/h \rfloor} & = \sum_{m=1}^{\lfloor (n-1)/h \rfloor}{\frac{b_m^2}{\mu^2}}  
\leq \frac{n-1}{h} \, \frac{(2Cn^{0.8190})^2}{\mu^2}  
\leq \frac{n}{n^{0.64}} \, \frac{4C^2n^{1.6380}}{(1/4)n^2} \leq 16C^2 n^{-0.0020} \cr
& \leq 2^{-9}c^6 \cr \endsplit$$
if $n=2^k$ is sufficiently large. So $S_n$ satisfies the assumptions of Lemma 3.5 with $c=1/4$ and 
$h:=n^{0.64}$ if $n=2^k$ is sufficiently large, indeed. Thus Lemma 3.5 implies that  
$${\Cal N}(S_n,v) > Ah^{-1}c^5n = Ac^5n^{0.36}$$
whenever $v$ is real with $|v| \leq 2^{-5}c^3 = 2^{-11}$ and $n=2^k$ is sufficiently large. 
\qed \enddemo

\demo{Proof of Theorem 2.8}
Suppose $P_k$ does not have a zero in the disk 
$$\left \{z \in {\Bbb C}: |z - e^{it_0}| < \frac{c_4}{n} \right \}\,.$$ 
Observe that
$$\split |e^{it} - e^{it_0}| & = |e^{it_0}(1 - e^{i(t-t_0})| = 
|t-t_0|\left|\sum_{j=1}^{\infty}{\frac{(i(t-t_0))^{j-1}}{j!}}\right| \cr
& \leq 2|t-t_0|\,, \qquad |t-t_0| \leq 1\,, \cr \endsplit$$
implies that $R_k \in {\Cal T}_n$ defined by
$R_k(t) = P_k(e^{it})P_k(e^{-it})$ does not have a zero in 
$$\left \{t \in {\Bbb C}: |t-t_0| < \frac{c_4/2}{n} \right\}\,.$$  
It follows from Lemma 3.6 and $\|R_k\|_K \leq 2n$ that
$$|R_k^{\prime}(t_0)| \leq 5e \sqrt{\frac{2n}{(c_4/2)/n}} \, \|R_k\|_K \leq 
\frac{20e}{\sqrt{c_4}} \, n^2 < cn^2$$ 
whenever
$$0 < c_4 < \frac{c^2}{400e^2}\,.$$
Hence, if we chose $c_4 > 0$ as above, $P_k$ must have a zero in the disk 
$$\left \{z \in {\Bbb C}: |z - e^{it_0}| < \frac{c_4}{n} \right \}$$   
whenever $R_k^{\prime}(t_0) \geq cn^2$.
\qed \enddemo

\demo{Proof of Theorem 2.7}
$$I_j := \left[ \frac{(2j-2)\pi}{n},\frac{2j\pi}{n} \right)\,, \qquad j=1,2,\ldots,n\,.$$
Let $\gamma := \sin^2(\pi/8)$ as before. By Saffari's Conjecture 1.1 proved by Rodgers [R-16] we have
$$m(\{t \in K: R_k(t) \leq \gamma n\}) > 2\pi(\gamma/4)$$
for every sufficiently large $n$. Hence, with the notation
$$A := \{t \in K: R_k(t) \leq \gamma n\}\,,$$
there are at least $n\gamma/4$ distinct values of $j \in \{1,2,\ldots,n\}$ such that
$A \cap I_j \neq \emptyset$ for every sufficiently large $n$.
On the other hand, by Lemma 3.2, for each $j \in \{1,2,\ldots,n\}$
there is a $t_j \in I_j$ such that $R_k(t_j) \geq 2\gamma n$. Hence by the Mean Value Theorem 
there are at least $n\gamma/4$ distinct values of $j \in \{1,2,\ldots,n\}$
for which there is a $\tau_j \in I_j$ such that 
$$R_k^{\prime}(\tau_j) \geq \gamma n (2\pi/n)^{-1} \geq \frac{\gamma}{2\pi}n^2$$
for every sufficiently large $n$. Hence, by Theorem 2.8, there are at least $n\gamma/4$ distinct values of 
$j \in \{1,2,\ldots,n\}$ such that the open disk $D_j$ centered at $e^{i\tau_j}$ of radius $c_4n^{-1}$ has 
at least one zero of $P_k$, where the absolute constant $c_4 > 0$ is chosen to $c := \gamma/(2\pi)$ as in 
the proof of Theorem 2.8, that is, $$0 < c_4 < \frac{\gamma^2}{1600\pi^2e^2}\,.$$   
\qed \enddemo

\head 5. Acknowledgement \endhead
The author thanks Stephen Choi for checking the details of the proof in this paper 
and for his computations leading to a correctly justified upper bound, stated 
as Lemma 3.4, for the autocorrelation coefficients of the Rudin-Shapiro polynomials.


\Refs \widestnumber\key{ABCD2}

\medskip

\ref \no AC-17
\by J.-P. Allouche, K.-K. S. Choi, A. Denise, T. Erd\'elyi, and B. Saffari
\paper Bounds on autocorrelation coefficients of Rudin-Shapiro polynomials
\paperinfo manuscript
\endref

\medskip

\ref \no BN-73 \by E. Beller and D.J. Newman,
\paper An extremal problem for the geometric mean of polynomials
\jour Proc. Amer. Math. Soc. \vol 39 \yr 1973 \pages 313--317
\endref

\medskip

\ref \no BP-32 \by A. Bloch and G. P\'olya,
\paper On the roots of certain algebraic equations
\jour Proc. London Math. Soc. \vol 33 \yr 1932 \pages 102--114
\endref

\medskip

\ref \no B-02 \by P. Borwein
\book Computational Excursions in Analysis and Number Theory
\publ Springer \publaddr New York \yr 2002
\endref



\medskip

\ref \no BE-95 \by  P. Borwein and T. Erd\'elyi
\book Polynomials and Polynomial Inequalities
\publ Springer \publaddr New York \yr 1995
\endref



\medskip

\ref \no BE-97 
\by P. Borwein and T. Erd\'elyi
\paper On the zeros of polynomials with restricted coefficients
\jour Illinois J. Math. \vol 41 \yr 1997 \pages no. 4, 667--675
\endref

\medskip

\ref \no BE-01
\by P. Borwein and T. Erd\'elyi
\paper Trigonometric polynomials with many real zeros and a Littlewood-type problem 
\jour Proc. Amer. Math. Soc. \vol 129 \yr no. 3, 2001 \pages 725--730
\endref

\medskip

\ref \no BE-07 
\by P. Borwein and T. Erd\'elyi
\paper Lower bounds for the number of zeros of cosine polynomials in the period: a problem of Littlewood
\jour Acta Arith. \vol 128 \yr 2007 \pages no. 4, 377--384
\endref

\medskip

\ref \no BE-08a 
\by P. Borwein, T. Erd\'elyi, R. Ferguson, and R. Lockhart
\paper On the zeros of cosine polynomials: solution to a problem of Littlewood
\jour Ann. Math. Ann. (2) \vol 167 \yr 2008 \pages no. 3, 1109--1117
\endref

\medskip

\ref \no BE-99 \by  P. Borwein, T. Erd\'elyi, and G. K\'os
\paper Littlewood-type problems on $[0,1]$
\jour Proc. London Math. Soc. \vol 79 \yr 1999 \pages 22--46
\endref

\medskip

\ref \no BE-13 \by  P. Borwein, T. Erd\'elyi, and G. K\'os
\paper The multiplicity of the zero at  $1$ of polynomials with constrained coefficients
\jour Acta Arith.  \vol 159 \yr 2013 \pages no. 4, 387--395
\endref

\medskip

\ref \no BE-08b 
\by P. Borwein, T. Erd\'elyi, and F. Littmann
\paper Zeros of polynomials with finitely many different coefficients
\jour Trans. Amer. Math. Soc. \vol 360 \yr 2008 \pages 5145--5154
\endref

\medskip

\ref \no BL-01 
\by P. Borwein and R. Lockhart
\paper The expected $L_p$ norm of random polynomials
\jour Proc. Amer. Math. Soc. \vol 129 \yr 2001 \pages 1463--1472
\endref

\medskip

\ref \no BM-00 \by P. Borwein and M.J. Mossinghoff
\paper Rudin-Shapiro like polynomials in $L_4$
\jour Math. Comp. \vol 69 \yr 2000 \pages 1157--1166
\endref

\medskip

\ref \no B-97
\by D. Boyd
\paper On a problem of Byrne's concerning polynomials with restricted coefficients
\jour Math. Comput. \vol 66  \yr 1997 \pages 1697--1703
\endref

\medskip

\ref \no B-73 
\by J. Brillhart
\paper On the Rudin-Shapiro polynomials 
\jour Duke Math. J. \vol 40 \yr 1973 \pages no. 2, 335--353
\endref

\medskip

\ref \no BL-76 \by J. Brillhart, J.S. Lemont, and P. Morton
\paper Cyclotomic properties of the Rudin-Shapiro polynomials
\jour J. Reine Angew. Math. (Crelle's J.) \vol 288 \yr 1976
\pages 37--65
\endref

\medskip

\ref \no CE-15a \by K.-K. S. Choi and T. Erd\'elyi
\paper Sums of monomials with large Mahler measure
\jour J. Approx. Theory \vol 197 \yr 2015 \pages 49--61
\endref

\medskip

\ref \no CE-15b \by K.-K. S. Choi and T. Erd\'elyi
\paper On a problem of Bourgain concerning the $L_p$ norms of exponential sums
\jour Math. Zeit. \vol 279 \yr 2015 \pages 577--584
\endref

\medskip

\ref \no CE-15c \by K.-K. S. Choi and T. Erd\'elyi
\paper On the average Mahler measures on Littlewood polynomials
\jour Proc. Amer. Math. Soc. Ser. B \vol 1 \yr 2015 \pages 105--120
\endref

\medskip

\ref \no CM-11 \by K.-K. S. Choi and M.J. Mossinghoff
\paper Average Mahler's measure and Lp norms of unimodular polynomials
\jour Pacific J. Math. \vol 252 \yr 2011 \pages no. 1, 31--50
\endref

\medskip

\ref\no CG-00 \by B. Conrey, A. Granville, B. Poonen, and K. Soundararajan
\paper Zeros of Fekete polynomials
\jour Ann. Inst. Fourier (Grenoble) \vol 50 \yr 2000 \pages 865--884
\endref

\medskip

\ref \no D-05 \by Ch. Doche
\paper Even moments of generalized Rudin-Shapiro polynomials
\jour Math. Comp. \vol 74 \yr 2005 \pages no. 252, 1923--1935
\endref

\medskip

\ref \no DH-04 \by Ch. Doche and L. Habsieger
\paper Moments of the Rudin-Shapiro polynomials
\jour J. Fourier Anal. Appl. \vol 10 \yr 2004 \pages no. 5, 497--505
\endref

\medskip

\ref \no D-08 \by P. Drungilas
\paper Unimodular roots of reciprocal Littlewood polynomials
\jour J. Korean Math. Soc. \vol 45 \yr 2008 \pages no. 3, 835--840
\endref

\medskip

\ref \no EZ-17 \by S.B. Ekhad and D. Zeilberger
\paper Integrals involving Rudin-Shapiro polynomials and sketch of a proof of Saffari's conjecture
\paperinfo to appear in the Proceedings of the Alladi60 conference
\endref

\medskip

\ref \no E-98 \by T. Erd\'elyi
\paper Markov-type inequalities for constrained polynomials with complex coefficients
\jour Illinois J. Math.\vol 42 \yr 1998 \pages 544--563
\endref

\medskip

\ref \no E-01 \by T. Erd\'elyi
\paper On the zeros of polynomials with Littlewood-type coefficient constraints
\jour Michigan Math. J. \vol 49 \yr 2001 \pages 97--111
\endref

\medskip

\ref \no E-08a \by T. Erd\'elyi
\paper An improvement of the Erd\H os-Tur\'an theorem on the distribution of zeros of polynomials
\jour C. R. Acad. Sci. Paris, Ser. I \vol 346 \yr 2008 \pages no. 5, 267--270
\endref

\medskip

\ref \no E-08b \by T. Erd\'elyi
\paper Extensions of the Bloch-P\'olya theorem on the number of real zeros of polynomials
\jour J. Th\'eor. Nombres Bordeaux \vol 20 \yr 2008 \pages no. 2, 281–-287
\endref

\medskip

\ref \no E-11 \by T. Erd\'elyi
\paper Sieve-type lower bounds for the Mahler measure of polynomials on subarcs
\jour Computational Methods and Function Theory \vol 11 \yr 2011 \pages 213--228
\endref

\medskip

\ref \no E-12 \by T. Erd\'elyi
\paper Upper bounds for the Lq norm of Fekete polynomials on subarcs
\jour Acta Arith.\vol 153 \yr 2012 \pages no. 1, 81--91
\endref

\medskip

\ref \no E-16a \by T. Erd\'elyi
\paper Coppersmith-Rivlin type inequalities and the order of vanishing of polynomials at $1$
\jour Acta Arith. \vol 172 \yr 2016 \pages no. 3, 271--284
\endref

\medskip

\ref \no E-16b \by T. Erd\'elyi
\paper On the number of unimodular zeros of self-reciprocal polynomials with coefficients from a finite set
\jour Acta Arith. \toappear
\endref

\medskip

\ref \no E-16c \by T. Erd\'elyi
\paper The Mahler measure of the Rudin-Shapiro polynomials
\jour Constr. Approx. \vol 43 \yr 2016 \pages no. 3, 357-369
\endref

\medskip

\ref \no E-17 \by T. Erd\'elyi
\paper Improved lower bound for the Mahler measure of the Fekete polynomials
\jour Constr. Approx. \toappear
\endref

\medskip

\ref \no E-18 \by T. Erd\'elyi
\paper The asymptotic value of the Mahler measure of the Rudin-Shapiro polynomials
\jour https://arxiv.org/abs/1708.01189 \paperinfo submitted
\endref

\medskip

\ref \no EL-07 \by T. Erd\'elyi and D. Lubinsky
\paper Large sieve inequalities via subharmonic methods and the Mahler measure of
Fekete polynomials \jour Canad. J. Math. \vol 59 \yr 2007 \pages 730--741
\endref

\medskip

\ref \no E-02 \by T. Erd\'elyi and J. Szabados
\paper Bernstein inequalities for polynomials with constrained roots
\jour Acta Sci. Math. (Szeged) \vol 68 \yr 2002 \pages no. 3-4, 937-952
\paperinfo corrected reprint of Acta Sci. Math. (Szeged) 68 (2002), no. 1-2, 163-178.
\endref

\medskip

\ref \no G-51 \by M.J. Golay
\paper Static multislit spectrometry and its application to the panoramic display of infrared spectra,
\jour J. Opt. Soc. America \vol 41 \yr 1951 \pages 468--472
\endref

\medskip

\ref \no HL-52 \by G.H. Hardy, J. E. Littlewood, and G. P\'olya
\book Inequalities \publ Cambridge Univ. Press
\publaddr London \yr 1952
\endref

\medskip

\ref \no L-61 \by J.E. Littlewood
\paper On the mean values of certain trigonometrical polynomials
\jour  J. London Math. Soc. \vol 36 \yr 1961 \pages 307--334
\endref

\medskip

\ref \no L-64 \by J.E. Littlewood
\paper On the real roots of real trigonometrical polynomials (II)
\jour  J. London Math. Soc. \vol 39  \yr 1964 \pages 511--552
\endref

\medskip

\ref \no L-66a \by J.E. Littlewood
\paper The real zeros and value distributions of real trigonometrical polynomials
\jour  J. London Math. Soc. \vol 41 \yr 1966 \pages 336-342
\endref

\medskip

\ref \no L-66b \by J.E. Littlewood
\paper On polynomials $\sum \pm z^m$ and $\sum e^{\alpha_{m}i} z^m$, $z=e^{\theta i}$
\jour  J. London Math. Soc. \vol 41 \yr 1966 \pages 367--376
\endref

\medskip

\ref \no L-68 \by J.E. Littlewood
\book Some Problems in Real and Complex Analysis
\publ Heath Mathematical Monographs \publaddr Lexington, Massachusetts \yr 1968
\endref

\medskip

\ref \no M-06a 
\by I.D. Mercer
\paper Unimodular roots of special Littlewood polynomials \jour Canad. Math. Bull.
\vol 49 \yr 2006 \pages no. 3, 438--447
\endref

\medskip

\ref \no M-80 \by H.L. Montgomery
\paper An exponential polynomial formed with the Legendre symbol
\jour Acta Arith. \vol 37 \yr 1980 \pages 375--380
\endref

\medskip

\ref \no M-06b \by K. Mukunda
\paper Littlewood Pisot numbers
\jour J. Number Theory \vol 117 \yr 2006 \pages no. 1, 106--121
\endref

\medskip

\ref \no R-16 \by B. Rodgers
\paper On the distribution of Rudin-Shapiro polynomials and lacunary walks on $SU(2)$
\jour arxiv.org/abs/1606.01637 \paperinfo to appear in Adv. Math.
\endref

\medskip

\ref \no Sch-32 
\by E. Schmidt
\paper \"Uber algebraische Gleichungen vom P\'olya-Bloch-Typos
\jour Sitz. Preuss. Akad. Wiss., Phys.-Math. Kl.  \yr 1932 \pages 321
\endref

\medskip

\ref \no Sch-33
\by I. Schur
\paper Untersuchungen \"uber algebraische Gleichungen
\jour Sitz. Preuss. Akad. Wiss., Phys.-Math. Kl. \yr 1933 \pages 403--428
\endref

\medskip

\ref \no S-51 \by H.S. Shapiro
\book Extremal problems for polynomials and power series
\publ Master thesis \publaddr MIT \yr 1951
\endref

\medskip

\ref \no Sz-34 
\by G. Szeg\H o
\paper Bemerkungen zu einem Satz von E. Schmidt uber algebraische Gleichungen
\jour Sitz. Preuss. Akad. Wiss., Phys.-Math. Kl. \yr 1934 \pages 86--98
\endref

\medskip

\ref \no T-96
\by M. Taghavi
\paper An estimate on the correlation coefficients of the Rudin-Shapiro polynomials 
\jour IJST \vol 20 \year 1996 \pages no. 2, Trans. A Sci., 235--240
\endref

\medskip

\ref \no T-97
\by M. Taghavi
\paper Upper bounds for the autocorrelation coefficients of the Rudin-Shapiro polynomials
\jour Korean J. Com. \& Appl. Math. \vol 4 \yr 1997 \pages no. 1, 39--46
\endref

\medskip

\ref \no TV-07
\by V. Totik and P. Varj\'u  \paper Polynomials with prescribed zeros and small norm 
\jour Acta Sci. Math. (Szeged) \vol 73 \yr 2007 \pages 593--612
\endref

\endRefs
\enddocument